\documentclass[12pt]{article} 
\usepackage{amsthm,amsmath,amssymb,enumerate,graphicx,pgf,tikz} 

\usepackage[all,2cell]{xy}          

\usepackage[utf8]{inputenc}
\usepackage{mathrsfs}
\usepackage{enumerate}
\usepackage{mathptmx}    
\newtheorem {coro}{Corollaire}

\newtheorem{lem}{Lemme}
\newtheorem{thm}{Théorème}

\newtheorem{exem}{Exemple}

\newtheorem{remarks}{Remarques}

\newtheorem{definition}{D\'efinition}

\newcommand{\ra}{\longrightarrow}

\newcommand{\st}{\stackrel}

\newcommand{\mm}{\mathscr{M}}
\newcommand{\md}{\mathscr{D}}

\newcommand{\N}{\mathbb N}

\newcommand{\p}{\mathscr{P}}
\newcommand{\mc}{\mathscr{C}}
\newcommand{\ms}{\mathscr{S}}

\newcommand{\ds}{\displaystyle}

\newcommand{\ot}{\otimes}

\newcommand{\map}{\mbox{map}}

\newcommand{\oppair}[4]{%
\xymatrix{%
#1 \ar@<.5ex>[r]^-{#3} &{#2}\ar@<.5ex>[l]^-{#4}%
}}
\newcommand{\oppairi}[4]{%
\xymatrix@1{%
#1 \ar@<.5ex>[r]^{#3} &{#2}\ar@<.5ex>[l]^{#4}%
}}
\DeclareRobustCommand{\coprod}{\mathop{\text{\fakecoprod}}}
\newcommand{\fakecoprod}{%
  \sbox0{$\prod$}%
  \smash{\raisebox{\dimexpr.9625\depth-\dp0}{\scalebox{1}[-1]{$\prod$}}}%
  \vphantom{$\prod$}%
}
  
\begin{document}
\title{Lemme de Yoneda pour les foncteurs \`a valeurs monoidales \footnote{Expos\'e aux séminaires de LATAO, Facult\'e des Sciences De Tunis, 8-3-2018}}
\author{Fethi Kadhi\footnote{\'Ecole Nationale des Sciences de l'Informatique, 	Manouba University.}       }    
  
\date{LATAO (FST): 8-03-2018}
\maketitle

\begin{abstract}
We consider a closed symmetric monoidal category $\mm$. We show that if $I$ is a small category then $\mm^I$ is a closed $\mm$-module.
We rewrite the Yoneda Lemma in the case of monoidal valued functors. We derive an adjoint functor theorem and we show that $\mm^I$ is a closed symmetric monoidal category
\end{abstract}

%


\tableofcontents

\section{Introduction}\label{intro}
On considère une catégorie monoïdale symétrique fermée $\mm$ et une petite catégorie $I$. On montre d'abord que la catégorie des foncteurs à valeurs dans $\mm$, notée $\mm^I$, possède une structure de $\mm$-module fermé. Cette structure nous permet, dans la suite, d'écrire le Lemme de Yoneda dans le cas des foncteurs à valeurs monoïdales et nous en déduisons un résultat de coreduction. Nous considérons un  foncteur de $\mm$-modules  $F:\mm^I\longrightarrow\mc$, nous démontrons que $F$ est un adjoint à gauche si et seulement si $F$ est cocontinu, dans un tel cas, nous donnons explicitement son adjoint à droite. Nous présentons quelques applications de ce théorème de foncteur adjoint. En particulier, nous démontrons que si $\mm$ est une catégorie monoïdale symétrique fermée alors la catégorie des foncteurs $\mm^I$ possède elle même une structure de catégorie monoïdale symétrique fermée.
\section{Catégories monoïdales}
La structure de catégorie monoïdale apparaît comme la catégorification d'un monoïde.
En informatique théorique le monoide libre est la structure la plus utilisée. Dans l'art abstrait le mot monoïde a été
 introduit par le peintre Jean Claude Bédard dans son livre: Pour un art schématique: Étude d'un monoïde graphique. En algèbre, un monoïde est la donnée d'un ensemble $E$ munie d'une loi de composition interne $*$ associative et admettant un élément neutre $e$. Formellement le triplet $(E,*,e)$ est un monoïde si
\begin{enumerate}
	\item $\forall x,y\in E, x*y\in E$.
	\item $\forall x,y,z\in E, (x*y)*z=x*(y*z)$.
	\item $\forall x\in E, x*e=e*x=x.$
	\end{enumerate}
	Les exemples classiques de monoïdes sont\\
	$(\N,+,0),(\N,\times,1), (\N,\max,0), (\p(E),\cup,\emptyset), (\p(E),\cap, E)\dots$\\
	Pour définir un monoïde dans le langage de la théorie des catégories on doit éviter au maximum le mot  "élément" et exprimer les axiomes par l'existence de flèches. Soit $A$ un objet de la catégorie $\ms$. On dit que $A$ est un monoïde lorsqu'il existe une flèche associatrice $\mu:A\times A\ra A$ et une flèche unitaire $\eta:1\ra A$ telles que les deux diagrammes suivants commutent:
	\begin{itemize}
	\item 
	\[
\xymatrix{
A\times A\times A \ar[rr]^{\mu\times 1_A} \ar[d]_{1_A\times \mu}&  &A\times A\ar[d]_{\mu}     \\
A\times A \ar[rr]^{\mu}&& A 
}
\]

 Ce diagramme signifie que la flèche $\mu$ induit une loi associative sur $A$.
\item \[
\xymatrix{
A\cong 1\times A \ar[r]^{\eta\times 1_A} \ar[rd]_{1_A}   &A\times A\ar[d]^{\mu} & A\times 1\cong A\ar[l]_{1_A\times\eta} \ar[ld]^{1_A}     \\
                        &A&
}
\]
Ce diagramme capte l'existence d'un élément neutre $e=\eta(1)$.
	
	\end{itemize}
	Maintenant nous pouvons abstraire cette d\'efinition pour découvrir la structure de catégorie monoïdale:\\
	On dit qu'une catégorie $\mm$ est monoïdale s'il existent
		
		\begin{enumerate}
			\item un foncteur
			$$ \begin{array}{rcl}
\otimes:\mm\times\mm&\ra & \mm\\
 (m,n)&\longmapsto&m\otimes n,										
\end{array}$$
\item un objet $k$ qui s'appelle objet unité,
 \item un isomorphisme  de foncteurs $a:(-\otimes-)\otimes-\ra-(\otimes-\otimes-)$ qui s'appelle associateur.
\item un isomorphisme  de foncteurs $r:-\otimes k\ra -$ qui s'appelle unitaire \`a droite,
\item un isomorphisme  de foncteurs $l:k\otimes -\ra -$ qui s'appelle unitaire \`a gauche
\end{enumerate}
tels que les diagrammes de cohérence suivants commutent:
 
\begin{itemize}

	\item \[
\xymatrix{
((m\otimes n)\otimes p)\otimes q \ar[r]^{a} \ar[d]_{a\otimes q}  & (m\otimes n)\otimes (p\otimes q)\ar[r]^{a}
&m\otimes(n\otimes (p\otimes q))     \\
(m\otimes(n\otimes p))\otimes q \ar[rr]_{a}&              &m\otimes((n\otimes p)\otimes q) \ar[u]_{m\otimes a}  
}
\]
$m,n,p,q\in\mm$. Ce diagramme s'appelle identit\'e de pentagone. 
\item \[
\xymatrix{
(m\otimes k)\otimes n \ar[rr]^{a} \ar[rd]_{r\otimes n}   &       &m\otimes(k\otimes n) \ar[ld]^{m\otimes l}     \\
                        &m\otimes c
}
\]
$m,n\in\mm$. Ce diagramme s'appelle identit\'e triangulaire.
\end{itemize}
L'\'ecriture $(\mm,\otimes,k)$ résume le fait que $\mm$ est une catégorie monoïdale pour le produit tensoriel $\otimes$ et l'objet unité $k$. L'exemple fondamental de catégorie monoïdale est la catégorie des ensembles $\ms$, dans ce cas le produit tensoriel est le produit cartésien, l'objet unité est le singleton, l'associateur $a$ est l'isomorphisme canonique de composante, suivant le triplet $(e,f,g)$, la fonction
$$ \begin{array}{cccl}
a:&(e\times f)\times g&\ra e\times(f\times g)\\
 &[(a,b),c]&\longmapsto[a,(b\times c)]										
\end{array}$$
 Un autre exemple typique de catégorie monoïdale est la catégorie $Vect_k$ des espaces vectoriels sur un corps $k$. Dans $Vect_k$, le produit tensoriel est donn\'e par le produit tensoriel d'espaces vectoriels, l'objet unité est le corps $k$ regardé comme un $k$-espace vectoriel. En théorie des catégories le produit tensoriel de deux espaces vectoriels $e$ et $f$ est l'objet de $Vect_k$ qui représente le foncteur covariant
 $$ \begin{array}{cccl}
B(e\times f,-):&Vect_k&\ra & \ms\\
 &w&\longmapsto&B(e\times f,w),										
\end{array}$$
avec $B(e\times f,w)$ est l'ensemble des applications bilinéaires de $e\times f$ dans $w$. Il est facile de vérifier que 
$B(e\times f,w)\cong Vect_k(e,Vect_k(f,w))$ et que $Vect_k(k,e)\cong e$.
 Ainsi, la propriété universelle qui identifie $e\otimes f$ s'\'ecrit:
$$Vect_k(E\otimes F,w)\cong B(E\times F,W)$$
Cette caractérisation du produit tensoriel de deux espaces vectoriels nous suffit, en utilisant les outils de la théorie des catégorie, pour vérifier l'existence d'un associateur $a$. En effet
Soient $e,f,g,w\in Vect_k$, on a
\begin{eqnarray*}
Vect_k((e\otimes f)\otimes g,w)&\cong& B((e\times f)\times g,w)\\
                                &\cong& B(e\times (f\times g),w)\\
																&\cong&Vect_k(e\otimes (f\otimes g),w)
	\end{eqnarray*}															
Par le principe de Yoneda, $(e\otimes f)\otimes g\cong e\otimes (f\otimes g)$ naturellement en $(e,f,g)$.
De plus 
\begin{eqnarray*}
Vect_k(e\otimes k,w)&\cong& B(e\times k,w)\\
                    &\cong& Vect_k[e, Vect_k(k, w)]\\
										&\cong&Vect_k(e,w)
	\end{eqnarray*}
Il s'ensuit, toujours par le même principe de Yoneda, qu'il existe un isomorphisme naturel, $r:e\otimes k\ra e$.
L'unitaire \`a gauche existe par le même argument.\\
Une catégorie monoïdale $(\mm,\otimes,k)$ est dite fermée s'il existe un internal-hom foncteur
$$ \begin{array}{cccl}
[-,-]:&\mm^{op}\times\mm&\ra & \mm\\
 &(a,b)&\longmapsto&[a,b]									
\end{array}$$
tel que $$\mm(a\otimes b,c)\cong\mm(a,[b,c])$$
naturellement en $a$ et $c$.
La cat\'egorie $\ms$ a une structure de catégorie monoïdale fermée. Le internal-hom foncteur est donné par
$[b,c]=\ms(b,c)=b^a$. Il est facile de vérifier que
$$\ms(a\times b,c)\cong \ms(a,c^b).$$
Dans la catégorie $Vect_k$, l'ensemble des applications linéaires de $e$ vers $f$ est lui même un espace vectoriel qui satisfait
l'isomorphisme $$Vect_k(e\otimes f,g)\cong Vect_k(e,Vect_k(f,g)).$$
ce qui rend $Vect_k$ une catégorie monoïdale fermée.\\
Une catégorie monoïdale $(\mm,\otimes,k)$ est dite symétrique s'il existe un isomorphisme naturel $s:m\otimes n\longmapsto n\otimes m$ qui est égal \`a son inverse et qui est compatible avec les diagramme de cohérence de la monoidalité.
D'une manière générale toute catégorie qui a une structure $ccc$ de closed cartesian category est une catégorie monoïdale symétrique fermée. En particulier, si $\mc$ est une catégorie localement petite alors la catégorie $\hat{\mc}$ des prefaisceaux de $\mc$ a une structure de catégorie monoïdale symétrique fermée. La catégorie $Cat$ des catégories localement petites a une structure de $ccc$ ce qui justifie l'utilisation de la notation exponentielle  $\mc^{\md}$ pour désigner la catégories de foncteurs de $\mc$ à $\md$.
\section{Structure de $\mm$- module}
Soit $(\mm,\otimes,k)$ une catégorie monoïdale symétrique fermée. $\mm$ est équipée  de trois isomorphismes naturels
$$ \begin{array}{rcl}
(m\otimes n)\otimes p& \st a \ra& m\otimes (n\otimes p)\\ 
k\otimes m & \st l \ra &m \\
m\otimes k&\st r \ra & m  						
\end{array}$$
\begin{definition}\label{do}
On dit qu'une catégorie $\mc$ a une structure de $\mm$-module s'il existe s'il existe un foncteur

$$ \begin{array}{lrcl}
\otimes&\mm\times\mc&\ra & \mc\\
 &(m,c)&\longmapsto&m\otimes c										
\end{array}$$
avec des isomorphismes naturels

 $$ (m\otimes n)\otimes c\st \alpha \ra m\otimes( n\otimes c)$$
et
$$k\otimes c\st \lambda \ra c$$
$m,n\in\mm$, $c\in\mc$
de sorte que les diagrammes de cohérence suivants commutent
\begin{itemize}
	\item \[
\xymatrix{
((m\otimes n)\otimes p)\otimes c \ar[r]^{\alpha} \ar[d]_{a\otimes c}  & (m\otimes n)\otimes (p\otimes c)\ar[r]^{\alpha}
&m\otimes(n\otimes (p\otimes c))     \\
(m\otimes(n\otimes p))\otimes c \ar[rr]_{\alpha}&              &m\otimes((n\otimes p)\otimes c) \ar[u]_{m\otimes\alpha}  
}
\]

$m,n,p\in\mm$ and $c\in\mc$
\item \[
\xymatrix{
(k\otimes m)\otimes c \ar[rr]^{\alpha} \ar[rd]_{l\otimes c}   &       &k\otimes(m\otimes c) \ar[ld]^{\lambda}     \\
                        &m\otimes c
}
\]
$m\in\mm$, $c\in\mc$
x
\item \[
\xymatrix{
(m\otimes k)\otimes c \ar[rr]^{\alpha} \ar[rd]_{r\otimes c}   &       &m\otimes(k\otimes c) \ar[ld]^{m\otimes\lambda}     \\
                        &m\otimes c
}
\]
$m\in\mm$, $c\in\mc$
\end{itemize}
On dit qu'un $\mm$-module $\mc$ est fermé s'il existe encore deux foncteurs:
$$\begin{array}{rrcl}
\map_{\mc}:&\mc^{op}\times\mc&\ra & \mm     \\
            &(m,n)&\mapsto      &\map_{\mc}(m,n)\\
\end{array}						             
$$
and
$$\begin{array}{rrcl}
\mm^{op}\times\mc&\ra & \mc    \\
            (m,d)&\mapsto      &d^m\\
\end{array}						             
$$
avec des isomorphismes naturels:
$$\mc(m\otimes c,d)\cong\mm(m,\map_{\mc}(c,d))\cong\mc(c,d^m).$$
\end{definition}
\begin{remarks}  Nous observons que:\\ 
\begin{enumerate}  
	\item Le foncteur $\otimes$ d\'efinit une action de $\mm$ sur $\mc$. Cette action est notée de la même manière que le produit tensoriel de deux objets de $\mm$. On distingue les deux multiplications par les deux objets multipli\'es. Si $m.n\in\mm$ alors $m\otimes n$ est le produit tensoriel de $\mm$. Si $m\in\mm$ et $c\in\mc$ alors $m\otimes c$ d\'esigne l'action de $m$ sur $c$.
	\item $\mm$ est supposée symétrique. Il s'ensuit que tout $\mm$-module \`a gauche ou \`a droite est en fait un $\mm$-module bilatéral. 
	\item Si on remplace $\mc$ par $\mc^{op}$ dans la d\'efinition précédente alors on fait apparaitre une action fermée de $\mm$ sur $\mc^{op}$, c'est l'action duale de $\mm$ sur $\mc^{op}$.
\end{enumerate}
\end{remarks}
La catégorie $\mm$ elle même peut être regardée comme un $\mm$-module fermé. Dans ce cas le produit tensoriel est confondu avec l'action fermé de $\mm$. Lorsque $I$ est une catégorie petite nous allons voir, dans la suite, que la catégorie des foncteur \`a valeurs dans $\mm$ a une structure de $\mm$-module fermé. La description de cette structure nécessite le rappel de quelques propriétés des concepts fin et cofin (end and coend).\\
Soit $F:\mc^{op}\times\mc\ra X$ un bifoncteur. $F$ est contravariant et covariant respectivement par rapport \`a la première variable et la deuxième variable. Si $c\stackrel{f}{\longrightarrow}c'$ est une flèche dans $\mc$ alors la flèche
 ${(f,1_{c'})}:(c,c'){\longrightarrow}(c',c')$ donne naissance \`a une flèche ${F(f,1_{c'})}:F(c,c') {\longrightarrow} F(c,c')$ notée aussi $F(f,c)$, la flèche $(1_c,f):(c,c)\ra (c.c')$ donne naissance \`a une flèche ${F(c,f)}:F(c,c) {\longrightarrow} F(c,c')$.
\begin{definition}
Soit $w$ un objet de $X$. on dit que $e:w\ra F$ est un coin  de $F$ s'il existe une famille $e_c:w\ra F(c,c)$ de flèches dans $X$ telles que pour toute flèche $c\stackrel{f}{\longrightarrow}c'$ le diagramme suivant commute:
\[
\xymatrix{
w\ar[r]^{e_c}\ar[d]_{e_{c'}}&F(c,c)\ar[d]^{F(c,f)}\\
F(c',c')\ar[r]_{F(f,c')}&F(c,c')
}
\]
\end{definition} 
Une {\it fin} de $F$ est un {\it coin} universel $\bar{e}:\bar{w}\ra F$:
\[
\xymatrix{
\forall w\ar[d]_{\exists !\bar{f}}\ar[rd]^{e}&\\
\bar{w}\ar[r]^{\bar{e}}&F
}
\]
Tout coin $e:w\ra F$ se factorise de façon unique $e=\bar{e}\bar{f}$.
 On  note $\bar{w}$ par $$\int_{c\in\mc}F(c,c).$$
Une cofin est un cocoin universel, on la note par
$$\int^{c\in\mc}F(c,c).$$
En catégorie enrichie on peut interpréter la fin comme une limite et la cofin comme une colimite.
Il s'ensuit qu'un foncteur continu préserve la fin et un foncteur cocontinu préserve la cofin. En particulier, si
$F:\mc^{op}\times\mc\ra\md$ et $d\in\md$
alors
\begin{enumerate}
	\item $$\md(d,\int_{c\in\mc}F(c,c))\cong\int_{c\in\mc}\md(d,F(c,c))$$
	\item $$\md(\int^{c\in\mc}F(c,c),d)\cong\int_{c\in\mc}\md(F(c,c),d)$$
\end{enumerate}
\begin{exem}
Soient $\mc$ et $\md$  deux catégories localement petites et $M,N:\mc\ra\md$ deux foncteurs alors
$$\int_{c\in\mc}\md(Mc,Nc)\cong \mc^{\md}(M,N)$$
En effet, soit 
$$\begin{array}{rccl}
F:&\mc^{op}\times\mc&\ra& \ms\\
  & (a,b)&\longmapsto&\md(Ma,Nb)
	\end{array}$$
	Il est facile de vérifier qu'un coin de $F$ correspond bijectivement \`a une transformation naturelle de $M$ \`a $N$.
\end{exem}
Considérons maintenant une catégorie monoïdale symétrique fermée $(\mm,\otimes,k)$. Le hom-foncteur $[-,-]$ qui fait de $\mm$ une cat\'egorie fermée  est noté exponentiellement de sorte que pour $m,n,p\in\mm$ il existe un isomorphisme naturel
$$\mm(m\otimes n,p)\cong\mm(m,p^n)$$

Soit $I$ une petite catégorie. De façon naturelle, on peut définir une action de $\mm$ sur la catégorie $\mm^I$ des foncteurs de $I$ \`a valeurs dans $\mm$. Cette action est donnée par le foncteur suivant:
$$\begin{array}{crclccl}
 \ot:&\mm\times \mm^I & \ra & \mm^I& & &\\
 &(m, M) & \longmapsto  & m \otimes M:&I&\ra&\mm\\
      &  &              &             &i&\longmapsto& m\otimes M_i\\
\end{array}$$
Dans ce qui suit on suppose que $\mm$ est bicomplete de sorte que toutes les limites et les colimites existent. On peut donc définir le foncteur suivant:
$$\begin{array}{crcllll}
 \mbox{map}_{\mm^I}:&(\mm^I)^{op}\times \mm^I & \ra & \mm \\
 &(M, N) & \longmapsto  & \mbox{map}_{\mm^I}(M,N)=\int_{i\in I} N_i^{M_i}\\ \end{array}$$
La structure monoïdale fermée de $\mm$ permet, en outre, de définir le foncteur suivant:
$$\begin{array}{crcllll}
 &\mm^{op}\times \mm^I & \ra & \mm^I& & &\\
 &(m, M) & \longmapsto  & M^m:&I&\ra&\mm\\
   &     &              &             &i&\longmapsto&M^m_i\\
\end{array}$$
\begin{lem}\label{l1}
La catégorie $\mm^I$  a une structure de  $\mm$-module fermé.
\end{lem}
\begin{proof}
Soit $m\in\mm$, $M,N\in\mm^I$, d'une part on a
\begin{eqnarray*}
\mm^I(m\otimes M,N)&\cong&\int_{i\in I}\mm(m\otimes M_i,N_i)\\
                    &\cong&\int_{i\in I}\mm(M_i\otimes m,N_i)\\
                   &\cong&\int_{i\in I}\mm(M_i,N^m_i)\\
									&\cong&\mm^I(M,N^m)\\
	\end{eqnarray*}
		D'autre part, on a
	\begin{eqnarray*}	
								\mm^I(m\otimes M,N)&\cong&\int_{i\in I}\mm(m\otimes M_i,N_i)\\
									&\cong&\int_{i\in I}\mm(m,N_i^{M_i})\\
									&\cong &\mm(m,\int_{i\in I}N_i^{M_i})\\
									&\cong&\mm(m,{map}_{\mm^I}(M,N))
	\end{eqnarray*}
	Il s'ensuit que nous avons des isomorphismes naturels
	 $$\mm^I(m\otimes M,N)\cong\mm^I(M,N^m)\cong\mm(m,{map}_{\mm^I}(M,N))$$
qui font de $\mm^I$ un $\mm$-module fermé.
	\end{proof}

\begin{definition}(Foncteurs de $\mm$-modules)\\
Soit $\mc, \md$ deux $\mm$-modules. On dit qu'un foncteur $\mc\st F \ra\md$ est un foncteur de $\mm$-modules s'ils existent des isomorphismes naturel $m\otimes F(c)\st \mu\ra F(m\otimes c)$ pour $m\in\mm, c\in\mc$
tels que les diagrammes suivants commutent 
 \begin{itemize}

\item \[
\xymatrix{(m\otimes n)\otimes F(c) \ar[rr]^{\mu}\ar[d]_{\alpha}&              &F((m\otimes n)\otimes c)\ar[d]^{F(\alpha)}
 \\
m\otimes(n\otimes F(c)) \ar[r]_{m\otimes\mu}   & m\otimes F(n\otimes c)\ar[r]_{\mu}
&F(m\otimes(n\otimes c))    
}
\]
$m,n\in\mm$ and $c\in\mc$

\item \[
\xymatrix{
k\otimes F(c) \ar[r]^{\lambda} \ar[d]_{\mu}   &  F(c)        \\
                        F(k\otimes c)\ar[ru]_{F(\lambda)}&
}
\]
$c\in\mc$

\end{itemize}
 \end{definition}
\begin{remarks}\label{rq}
Soient $\mc$,  $\md$ deux $\mm$-modules fermés et $\oppairi{\mc}{\md}{F}{G}$ une adjonction. Il s'ensuit que
\begin{enumerate}
	\item $F$ est un foncteur de $\mm$-modules si et seulement si $G^{op}:\md^{op}\ra\mc^{op}$ est un foncteur de $\mm$-modules pour l'action duale de $\mm$ sur $\mc^{op}$ et $\md^{op}$, dans un tel cas l'adjonction $(F,G)$ s'appelle une adjonction de $\mm$-modules.
	
	\item\label{r0} Supposons que $(F,G)$ est une adjonction de $\mm$-modules alors $$\map_{\md}(F(c),d) \cong\map_{\mc}(c,G(d)),$$ naturellement en $c$ et $d$.\\
		En effet,	soient $c\in\mc$ et $d\in\md$.	Pour $m\in\mm$
	\begin{eqnarray*}
	\mm(m,\map_{\md}(F(c),d))&\cong&\md(m\otimes F(c),d)\\
	                          &\cong& \md(F(m\otimes c),d)\\
														&\cong& \mc(m\otimes c, G(d))\\
														&\cong&\mm(m,\map_{\mc}(c,G(d))
	\end{eqnarray*}
	Le résultat en découle par le principe de Yoneda. 
\end{enumerate}
\end{remarks}
%
%
%
%
%
%
%
%
%
%
%
%
\section{Lemme de Yoneda pour les catégories monoïdales}
Dans le cas de la catégories des foncteurs à valeurs dans $\ms$ le Lemme de Yoneda s'écrit
$$\ms^I(Hi,M)\cong Mi$$
avec $Hi(j)=I(i,j)$. Pour découvrir comment s'écrit le lemme de Yoneda dans le cas de foncteurs à valeurs monoïdales, nous avons besoin de trouver un pont entre $\ms^I$ et $\mm^I$. Cherchons d'abord une adjonction entre $\ms$ et $\mm$.
Soient $U:\ms\ra \mm$ et $V:\mm\ra \ms$ les foncteurs définis par $U(s)=\ds\coprod_{s}k$ et $V(m)=\mm(k,m).$
$\oppairi{\ms}{\mm}{U}{V}$ est une adjonction.\\
 En effet, pour $s\in\ms$ et $m\in\mm$ on a
	\begin{eqnarray*}
	\mm(U(s),m)&\cong&\mm(\coprod_sk,m)\\
	            &\cong& \prod_s\mm(k,m)\\
							&\cong& \mm(k,m)^s \\
							&\cong&\ms(s,\mm(k,m))\\
							&\cong&\ms(s,V(m))
	\end{eqnarray*}
	
	Cette adjonction entre $\ms$ et $\mm$ induit une adjonction entre $\ms^I$ et $\mm^I$:
	$$\oppairi{\ms^I}{\mm^I}{U^I}{V^I}$$
	Pour $F\in\ms^I$ et $G\in\mm^I$ on a $U^I(F)i=U(Fi)=\ds\coprod_{Fi}k$ et $V^I(G)i=\mm(k,Gi)$. On a
	\begin{eqnarray*}
\mm^I(U^I(F),G)&\cong&\int_{i\in I}\mm(U(Fi),Gi)\\
             &\cong&\int_{i\in I}\ms(Fi,V(Gi))\\
						&\cong&\ms^I(F,V^I(G))   \\
	\end{eqnarray*}
	Soient $i\in I$ et 	$h_i=U^I(Hi)=U \circ Hi $, Il s'ensuit
$$ \begin{array}{lll}
 h_i:&I \ra & \mm\\
 &j\longmapsto  & \ds\coprod_{I(i,j)}k																	
\end{array}$$ 
\begin{lem}\label{eq1}
Pour $M\in\mm^I$, on a $$\mm^I(h_i,M)\cong\mm(k,M_i)$$
\end{lem}
\begin{proof}
\begin{eqnarray*}
\mm^I(h_i,M)&\cong&\mm^I(U^I(Hi),M)\\
             &\cong&\ms^I(Hi,V^I(M)) \\
						&\cong&V^I(M)_i  \:\: \mbox{  Par le lemme de Yoneda dans $\ms$ }\\
						&\cong&V(M_i)\\
	          &\cong&\mm(k,M_i) 
\end{eqnarray*}
\end{proof}
Si $\mm$ est la catégorie des ensembles alors l'objet unité est le singleton noté 1. Le Lemme \ref{eq1} s'écrit alors
$\ms^I(h_i,M)\cong\ms(1,M_i)\cong M_i$, ainsi, on retrouve le lemme de Yoneda dans le cas classique.
Ce pendant  ce lemme donne un isomorphisme entre deux objets de $\ms$ et par suite on ne peut pas le considérer comme l'homologue du lemme de Yoneda dans le cas des foncteurs \`a valeurs monoïdales. Nous cherchons un lemme qui donne un isomorphisme entre deux objets de $\mm$. Nous disposons du foncteur $\mbox{map}_{\mm^I}$ qui a permis de montrer que $\mm^I$ a une structure de $\mm$-mdule fermé c'est grace \`a ce foncteur et au Lemme \ref{eq1} que nous \'ecrivons le lemme de Yoneda pour $\mm^I$.
\begin{lem}\label{l2}{(Lemme de Yoneda pour les catégories monoïdales symétriques fermées)}\\
$h_i$ \'etant comme dej\`a d\'efini. Soit $M\in\mm^I$, on a
$$\mbox{map}_{\mm^I}(h_i,M)\cong Mi$$
\end{lem}
\begin{proof} Pour $m\in\mm$, 
\begin{eqnarray*}
\mm(m, map_{\mm^I}(h_i,M))&\cong& \mm(m\otimes h_i,M)\\
                  				&\cong& \mm^I(h_i,M^m)\\
								         	&\cong& \mm(k,M_i^m) \mbox{ (par le  Lemme \ref{eq1}})\\
							        		&\cong&\mm(m,M_i)\\
		        	\end{eqnarray*} 
Tous les isomorphismes sont naturels, le résultat en découle par le principe de Yoneda. 
\end{proof}
Soit $i\in I$, on désigne par $Ev_i$ le foncteur  i-evalution defini par:
$$ \begin{array}{cccl}
Ev_i:&\mm^I&\ra & \mm\\
 &M&\longmapsto&Mi										
\end{array}$$
On a le resultat suivant:
\begin{lem}\label{nl3}
$Ev_i$ est un adjoint \`a droite, son adjoint \`a gauche est le foncteur
 $\mm \st {F_{i}}\ra \mm^{I}$ défini par $F_{i}(m)=h_{i} \ot m$.
\end{lem}
\begin{proof}
Il s'agit d'une application du Lemme \ref{l2}.
Soit $m \in \mm $ et $M \in \mm^{I}$, on a
\begin{eqnarray*}
 \mm ^{I}(F_{i}(m),M)&\cong &\mm^{I}(h_{i} \ot m ,M)\\
                     &\cong& \mm(m,map_{\mm^{I}}(h_{i},M))\\
										&\cong& \mm(m,Ev_{i}(M))\mbox{(par le  lemme \ref{l2})} 
\end{eqnarray*} 
 \end{proof}
\begin{definition}
 Soit $M: I\ra \mm$, on définit le codifferentiel de $M$ par 
$$ \begin{array}{rrcllll}
CM:&I^{op}\times I&\ra & \mm^I&&&  \\
 &(i,j)&\longmapsto & h_i\otimes M_j	:& I& \ra &\mm\\ 
&&&& q& \longmapsto &  h_i(q)\otimes M_j 												
\end{array}$$
\end{definition}
Le lemme suivant parait comme l'analogue, en théorie des catégories, du théorème fondamental du calcul. Cette analogie justifie le choix de notre terminologie
\begin{lem}\label{l0} 
Soit $ M\in \mm^{I}$,   $M$ est isomorphe au cofin (coend) de son codifferentiel  $CM$, i.e.  $$M \cong \int^{i\in I} CM(i,i)$$ 
\end{lem}
\begin{proof}
Il s'agit d'une autre application du Lemme \ref{l2}
 Pour $N \in \mm^{I}$, on a
\begin{eqnarray*}
\mm^I (\int^{i\in I} CM(i,i),N)&\cong&\int_{i\in I} \mm^I (CM(i,i),N)\\
                    &\cong&  \int_{i\in I} \mm^I (h_i\otimes M_i,N)\\
										&\cong&\int_{i\in I} \mm  ( M_i, map_{\mm^{I}} (h_i,N))\\
										&\cong&\int_{i\in I} \mm  ( M_i, N_i) (\mbox{ D'après le lemme \ref{l2}})\\ 
										&\cong&\mm^I (M,N)
\end{eqnarray*}
Il s'ensuit que  $M \cong \int^{i\in I} CM(i,i)$  
\end{proof}
Maintenant nous avons bien préparé le terrain pour énoncer un théorème de foncteur adjoint
\section{Théorème de foncteur adjoint}
Dans son livre \cite{ML2}, S. Mac Lane dit que pendant une douzaine d'années on a utilisé les catégories juste comme un langage, dans la suite cette vision a changé grâce à l'apparition des foncteurs adjoints. Son élève S. Awoday explique dans \cite{SA}
que la théorie des catégories est la seule discipline mathématique qui capte la notion d'adjonction. Généralement, un théorème de foncteur adjoint  fournit une condition qui garantit l'existence d'un adjoint pour un foncteur donné \cite{tl}. Le théorème suivant s'intéresse à un foncteur de $\mm$-modules.
\begin{thm}\label{nt1}
Soit $\mm^{I} \st {F}  \ra \mc $ un foncteur de $\mm$-modules. $F$ est un adjoint \`a gauche si et seulement si $F$ est cocontinu, dans un tel cas, son adjoint \`a droite $G$ est d\'efini par
$$G(Y)_{i}=map_{\mc}(F(h_{i}),Y), i\in I$$ 
\end{thm}
\begin{proof}
Un adjoint à gauche est nécessairement cocontinu. Réciproquement, supposons que $F$ est cocontinu. Si $F$ admet un adjoint à droite $G$ alors que devrait être $G$?
Soient $X \in \mm^I$ et $Y \in \mc$, d'après la remarque \ref{rq} .\ref{r0}
$$map_\mc(F(X), Y) \cong map_{\mm^I} (X, G(Y))$$
Prenons $X=h_i$ et appliquons le Lemme \ref{l2}, nous obtenons
 $$ G(Y)_i \cong map_\mc(F(h_i),Y)$$ 
Vérifions maintenant que l'objet fourni par cette formule est l'adjoint à droite de $F$. On a
\begin{eqnarray*}
\mm^{I}(X, G(Y))&\cong&\int_i \mm (X_i, G(Y)_i)\\
                    &\cong&\int_i  \mm  (X_i, map_ \mc (F(h_i), Y ))\\
                   &\cong&\int_i \mc (F(h_i)\otimes X_i,Y)\\
									&\cong&\int_i \mc (F(h_i \otimes X_i),Y)  \mbox{($F$ est un foncteur de $\mm $-modules)   }\\
									&\cong&\mc (\int^i F(h_i \otimes X_i),Y)\\
									&\cong&\mc (F(\int^i  h_i\otimes X_i),Y)  \mbox{(cocontinuité  de $F$)}\\
									&\cong&\mc (F(X),Y)  \mbox{(Lemma \ref{l0})} 
\end{eqnarray*}
 Ainsi $G$ est l'adjoint à droite de $F$.
\end{proof}

\begin{coro}\label{ex1}
Soit $ \st \Phi \ra J$ un foncteur de deux petites catégories $I$ et $J$, $\Phi$ induit un foncteur de $\mm$-modules $\mm^J \st F \ra \mm^I$. $F$ admet un adjoint \`a drote $G$ d\'efini par

$$G(Y)_j =   map_{\mm^{I}} (h_j\circ \phi,Y).$$ 
	\end{coro}
	\begin{proof}
	$F=-\circ \Phi$ est un foncteur de $\mm$-modules cocontinu. On obtient $G$ par application directe du théorème \cite{nt1} 
	\end{proof}
\begin{coro}(Structure monoïdale de  $\mm^I$)\\
Si $I$ est une categorie petite et $\mm$ est une catégorie bicomplete qui a une structure de catégorie monoidiale symétrique fermée alors $\mm^{I}$ a une structure de catégorie monoidiale symétrique fermée:  
\begin{enumerate}
	\item Le produit tensoriel est défini par:
	$$\begin{array}{rrcllll}
 \otimes:&\mm^I\times \mm^I & \ra & \mm^I& & &\\
         &(M, N) & \longmapsto  & M\otimes N:&I&\ra&\mm\\
         &      &              &             &i&\longmapsto& M_i\otimes N_i\\
	\end{array}$$
	\item L'objet unité est le foncteur
	$$\begin{array}{lll}
 K:&I \ra & \mm\\
   &i \longmapsto  & k
        	\end{array}$$
	\item Le hom-internal foncteur est d\'efini par:
	$$\begin{array}{ccllll}
(\mm^I)^{op}\times \mm^I & \ra & \mm^I& & &\\
        (N, P) & \longmapsto  & P^N:&I&\ra&\mm\\
              &              &             &i&\longmapsto& map_{\mm^{I}} (h_{i} \ot N,P)\\
\end{array}$$
\end{enumerate}
\end{coro}
\begin{proof}
	Pour $N\in \mm^{I}$ On d\'efinit le foncteur  $$F_{N}:\mm^{I} \ra \mm ^{I}
$$
                                                       $$M\longmapsto M \ot N$$
$F_{N} $ est un foncteur cocontinu de $\mm$-modules, par le théorème \ref{nt1}, $F_{N}$ admet un adjoint \`a droite $G_N$ donn\'e par $G_N(P)_{i}=map_{\mm^{I}} (h_{i} \ot N,P)$.

Ceci termine la preuve. Ce r\'edultat peut \^etre obtenu en utilisant la notion de {\it comonadic functor} et les extentions de Kan, voir \cite{nlab}. Dans ce corollaire, on remarque qu'on a pu obtenir le m\^eme r\'esultat avec des outils beaucoup plus simples.
\end{proof}

\end{document}